\documentclass[12pt]{article}
\usepackage[cp1251]{inputenc}
\usepackage{amssymb,amsmath}
\usepackage[T2A]{fontenc}
\begin{document}
\title{On modules over group rings of groups \\ 
with restrictions on  the system of all proper subgroups}
\author{O.Yu.Dashkova}
\date{\it \small Department of Mathematics and Mechanics, Dnepropetrovsk National University, 
prospekt Gagarina, 72, 49010, Dnepropetrovsk, Ukraine,
odashkova@yandex.ru}
\maketitle

\begin{abstract}
We consider the class   $\mathfrak M$  of $\bf R$--modules where $\bf R$ is an associative ring. Let  $A$ be a module over  a group ring  $\bf R$$G$ where $G$ is a group and  let
$\mathfrak L(G)$ be a set of all proper subgroups of $G$ such that if $H \in  \mathfrak L(G)$ then $A/C_{A}(H)$ belongs to  $\mathfrak M$.
We study an $\bf R$$G$--module $A$ such that $G \not = G'$, $C_{G}(A) = 1$,   $A/C_{A}(G) \not \in \mathfrak M$,
and $\mathfrak M$ is one of the classes: artinian $\bf R$--modules, 
 minimax $\bf R$--modules,   finite $\bf R$--modules. We consider  the cases: 1) $\mathfrak M$ is a class of all artinian $\bf R$--modules, $\bf R$ is either a ring  of integers or a ring of $p$--adic integers;
2) $\mathfrak M$ is a class of all minimax  $\bf R$--modules,  $\bf R$ is a ring  of integers, G is a locally soluble group; 3) $\mathfrak M$ is a class of all finite   $\bf R$--modules,  $\bf R$ is an associative ring.
In these   cases we prove that $G$ is isomorphic to a quasi--cyclic $q$--group  for some prime $q$.

\end{abstract}

$$ \textbf{Keywords: Module over group ring; artinian module; minimax module} $$

$$ \textbf{ AMS Subject Classification: 20F19, 20H25} $$

$$ \textbf{\S1. Introduction} $$

Investigation of modules over group rings is an important direction in  algebra. Modules over group rings of finite groups have been   studied 
by many authors. If $G$ is an infinite group, the situation is totally different. The study of this case requires some additional restrictions. 
Modules over group rings of infinite groups have been considered recently in \cite{KOS}. Artinian and noetherian modules over 
group rings are a broad class of modules over group rings.  Recall that a module is called artinian if the partially ordered set of all its submodules satisfies the minimal condition. A module is called
noetherian if the partially ordered set of all its submodules satisfies the maximal condition.
 Natural generalization of the classes of artinian and noetherian
modules is the class of minimax modules (chapter 7 \cite {KSS}). Let $A$ be an $\bf
R$--module,  $\bf R$ be  an associative ring.   $\bf
R$--module   $A$ is called minimax if it has the finite series of
submodules such that every its factor is either a noetherian $\bf R$--module   or
 an artinian $\bf R$--module. 
 It arises the question on investigation of  modules over group rings which are not artinian or  noetherian or else 
minimax but are similar to these modules in some sence. 
 
Let  $\mathfrak M$ be a class of $\bf R$--modules where $\bf R$ is an associative ring and let $A$ be a  module over a group ring  $\bf R$$G$ where $G$ is a group. Let
$\mathfrak L(G)$ be a set of subgroups of $G$ such that if $H \in $$\mathfrak L(G)$ then $A/C_{A}(H)$ belongs to  $\mathfrak M$.
  B.A.F. Wehrfritz have considered  groups $G$ of automorphisms of a module  $A$ over a ring  $\bf R$
 if  $\mathfrak M$ is  one of the  classes  of noetherian,  artinian  or  finite $\bf R$--modules and  $\mathfrak L(G)$ contains  all finitely generated subgroups of $G$  \cite {WBW}--\cite {FFW}, \cite {WWF}.  

Let $A$  be an  $\bf R$$G$--module   such that all proper subgroups of $G$ belong to $\mathfrak L(G)$ but $G$ does not belong to $\mathfrak L(G)$. 
If  $A$ is an $\bf R$$G$--module, $\bf R$$=F$  is a field of prime characteristic, $C_{G}(A) = 1$,  $G$ is an almost  locally soluble group  then $G$ is  isomorphic to a quasi--cyclic $q$--group  for some prime $q$ 
\cite {DEK}. It was investigated the case where  $A$ is an $\bf R$$G$--module, $\bf R$ is a ring $\mathbb  Z_{p^{\infty}}$ of $p$--adic integers, $C_{G}(A) = 1$,  $G$ is an infinite   soluble group and $\mathfrak M$  
is a class of all artinian $\bf R$--modules  \cite {ODD}.  It was proved that $G$  is  isomorphic to a quasi--cyclic $q$--group  for some prime $q$ also. In \cite {OYD} it was considered the case where $A$ is an $\bf R$$G$--module, $\bf R$ is 
a ring $\mathbb  Z$ of  integers,  $C_{G}(A) = 1$, $G$ is an infinite   soluble group and $\mathfrak M$  is a class of all artinian  $\bf R$--modules.
In this case $G$ is  isomorphic to a quasi--cyclic $q$--group  for some prime $q$ too.

In \cite{KSC} the authors   investigated the case where  $A$ is an $\bf R$$G$--module, $\bf R$ is 
a ring $\mathbb  Z$ of  integers, $C_{G}(A) = 1$,  $G$ is a locally generalized radical  group and  $\mathfrak M$  
is a class of all  artinian-by-(finite rank)  $\bf R$--modules.

We study $\bf R$$G$--module $A$ such that $G$ is an infinite group, $G \not = G'$, $C_{G}(A) = 1$, $A/C_{A}(G) \not \in  \mathfrak M$  and $\mathfrak M$ is one of the classes:  1) artinian $\bf R$--modules;
2) minimax $\bf R$--modules; 3) finite $\bf R$--modules. It should be noted that the class of groups which are different from their derived  subgroups  is sufficiently  broad.
 It contains all soluble groups, $ZD$--groups and free groups.
The main results of this work are theorems 2.1--2.3.

$$ \textbf{\S2. On modules over group rings of groups}$$
$$ \textbf{ with restrictions on the system of all proper subgroups } $$

Later we consider $\bf R$$G$--module $A$ such that $C_{G}(A) =1$, $ \bf R$ is an associative ring.  At first we prove some preliminary results.

{\bf Lemma 2.1.}
 Let $A$ be an $\bf R$$G$--module, $K$, $L$ be subgroups of $G$. If $A/C_{A}(K)$ and  $A/C_{A}(L)$ are minimax $\bf R$--modules then $A/C_{A}(\langle K, L \rangle )$  is 
 a minimax  $\bf R$--module also. 

{\bf Lemma 2.2.} 
 Let $A$ be an $\bf R$$G$--module where  $G$ is an infinite group, $G \not = G'$,  $ \bf R$ is an associative ring.  If $A/C_{A}(G)$ is not a minimax  
$\bf R$--module and $A/C_{A}(H)$ is a minimax $\bf R$--module for every proper subgroup $H$ of $G$ then  $G$ has not proper subgroups of finite index and $G/G'$ is isomorphic to a   quasi--cyclic $q$--group for some prime $q$.

{\bf Proof.}
 We prove that   $G$ is an infinite generated group.  Otherwise let  
$\{ x_{1}, x_{2}, \cdots, x_{m} \}$  be  a minimal system of generatings of  $G$. If  $m=1$ then $G$ is an infinite cyclic group. Therefore  $G$ is generated by two proper subgroups. By lemma 2.1  
$A/C_{A}(G)$ is a minimax  $\bf R$--module. Contradiction.  If  $k > 1$ then 
 $G$ is generated by proper subgroups $\langle x_{1}, x_{2}, \cdots, x_{m-1} \rangle$ and  $\langle x_{m} \rangle$. We have a contradiction also.  It follows that 
$G$ is an infinite generated group. Now  we prove 
 that $G$ has not proper subgroups of finite index. Otherwise if  $N$ is a proper subgroup of $G$ of  finite index   then we can choose a finitely generated subgroup 
$M$ such that  $G = MN$ where $M$ and  $N$ are proper subgroups of $G$. By lemma 2.1 $A/C_{A}(G)$ is  a minimax $\bf R$--module. Contradiction. 

 Let  $D$ be a derived subgroup of $G$. As  $G$ has not proper subgroups of finite index then $G/D$  is infinite. By lemma 2.1 an abelian quotient group  $G/D$ can not be generated by two proper 
subgroups. Let   $G/D$ be a  nonperiodic group and  $T/D$ be  a periodic part of $G/D$. Then  
$G/T$ is  generated by two proper subgroups.   Contradiction with lemma 2.1. 
 Therefore $G/D$ is periodic.  Hence  $G/D$ is isomorphic to a quasi--cyclic $q$--group for some prime $q$ (p.152 \cite {KUR}). The lemma is proved.

Lemmas 2.1 and  2.2 are  valid if the minimax condition is replaced by the artinian condition or the finiteness condition.

{\bf Lemma 2.3.}
Let $A$ be an artinian $\bf R$--module where $\bf R$$= \mathbb  Z_{p^{\infty}}$ is a ring of $p$--adic integers. Then the additive group of   $A$ is Chernikov and its divisible part is a  $p$--group.

{\bf Proof.}
 Let   $P$ be  a maximal ideal of  $\bf R$.  Then the  additive group of 
$\bf R$$/P$ has an order  $p$ and the additive group of 
$\bf R$$/P^{k}$ is a cyclic group of an orger 
$p^{k}$. Let    $\bf R $$/P^{k} = \langle a_{k}\rangle $,  $k=1, 2, \cdots$, 
 $\pi_{k}^{k+1}:\hspace{0,2cm} $$\bf R$$/P^{k} \longrightarrow $$\bf R$$/P^{k+1}$ where 
 $$
\pi_{k}^{k+1}(a_{k}) = pa_{k+1}
$$
 and  let $\pi_{k}^{m}: \hspace{0,2cm} $$\bf R$$/P^{k} \longrightarrow $$\bf R$$/P^{m}$ where 
 $m>k$ and 
$$
\pi_{k}^{m}(a_{k}) = p^{m-k}a_{m}.
$$ 
We consider  an injective limit of the set of 
 $\bf R$$/P^{k}$, $ k=1,2,\cdots ,n, \cdots \hspace{0,2cm}$. From the choice of  $a_{1}$  it follows that  $pa_{1} = 0$.
Therefore this   injective limit is isomorphic to a quasi--cyclic $p$--group  $C_{p^{\infty}}$. It follows that an additive group of a Pr$\ddot u$fer 
$\bf R$--module is isomorphic to a quasi--cyclic $p$--group  $C_{p^{\infty}}$ (ch. 5 \cite{KSS}).   By theorem 7.13 \cite{KSS} 
an artinian $\bf R $--module is decomposed in a direct sum  $A = a_{1} $$\bf R$$ \oplus a_{2} $$\bf R$$ \oplus \cdots \oplus a_{n} $$\bf R$$ \oplus C_{1} \oplus
 \cdots  \oplus C_{k}$ where  $C_{i}$ is a Pr$\ddot u$fer $P_{i}$--module,  $P_{i} \in  Spec($$\bf R$$), i=1, \cdots , k, Ann_{\bf R}(a_{j})=P^{m_{j}}_{j},
P_{j}  \in  Spec($$\bf R$$),  j=1, \cdots, n$.
 Every ideal of  $\bf R$  has a finite index in $\bf R$ (ch.  6 \cite {KUR}). Therefore $a_{j} $$\bf R$ is a finite $\bf R$--module  for each  $j=1, \cdots, n$. It follows that  the additive group of an artinian 
 $\mathbb  Z_{p^{\infty}}$--module  $A$ is Chernikov and its divisible part is a  $p$--group. 
The lemma is proved.

Now we prove the main results of this work.

{\bf Theorem 2.1.} 
 Let $A$ be an $\bf R$$G$--module where  $G$ is an infinite group, $G \not = G'$,  $\bf R$ is either a ring $ \mathbb  Z$ of integers or a ring $\mathbb  Z_{p^{\infty}}$ of $p$--adic integers.
If $A/C_{A}(G)$ is not an artinian $\bf R$--module and $A/C_{A}(H)$ is an artinian $\bf R$--module  for every proper subgroup $H$ of $G$ then $G$ is isomorphic to a quasi--cyclic $q$--group  $C_{q^{\infty}}$ for some prime $q$.

{\bf Proof.}
  Let  $D$ be a derived subgroup of $G$. By lemma 2.2    $G/D$ is isomorphic to a quasi--cyclic $q$--group for some prime $q$. Let $H/D$ be any finite subgroup of  $G/D$. 
Since  $H$ is a proper subgroup of  $G$ then $A/C_{A}(H)$ is an artinian $\bf R$--module. If $\bf R$ is a ring $ \mathbb  Z$ of integers then   $A/C_{A}(H)$ is an abelian group with the minimal condition for  subgroups. 
Therefore  $A/C_{A}(H)$ is a Chernikov group. If $\bf R$ is  a ring $\mathbb  Z_{p^{\infty}}$ of $p$--adic integers then $A/C_{A}(H)$ is a Chernikov group by lemma 2.3.
 It  follows that    $A/C_{A}(H)$ is the union of finite characteristic subgroups $A_{n}/C_{A}(H)$, $ n=1, 2, \cdots$,  and for each  $ n=1, 2, \cdots$, we have that 
$G/C_{G}(A_{n}/C_{A}(H))$ is  finite. By lemma 2.2  $G$ has not proper subgroups of finite index. Then  $G = C_{G}(A_{n}/C_{A}(H))$ for each $ n=1, 2, \cdots$. It follows that $[G, A_{n}] \leq C_{A}(H)$ for each $ n=1, 2, \cdots$. Therefore 
$[G,A] \leq C_{A}(H)$. From the choice of  $H$ it  follows that $[G,A] \leq C_{A}(G)$ and so  
$G$ acts trivially at every factor of the series 
$0 \leq C_{A}(G) \leq A$. By Kaluzhnin theorem (p. 144 \cite {KAM}) 
 $G$ is an abelian group. It  follows that   $G$ is isomorphic to a quasi--cyclic $q$--group  $C_{q^{\infty}}$ for some prime $q$.  The theorem is proved.

{\bf Theorem 2.2.} 
Let $A$ be a  $ \mathbb  Z$$G$--module where  $G$ is an infinite locally soluble  group, $G \not = G'$,  $ \mathbb  Z$ is a ring of  integers.
If $A/C_{A}(G)$ is not a minimax  $ \mathbb  Z$--module and $A/C_{A}(H)$ is a minimax $ \mathbb  Z$--module for every proper subgroup $H$ of $G$ then $G$ is isomorphic to a quasi--cyclic $q$--group  $C_{q^{\infty}}$ for some prime $q$.

{\bf Proof.}
  Let  $D$ be a derived subgroup of $G$. By lemma 2.2    $G/D$ is isomorphic to a quasi--cyclic $q$--group for some prime $q$. 
At first we consider the case where there is the proper subgroup  $L$ such that  $D \leq L$ and  $ A/C_{A}(L)$ is not an artinian $\mathbb{Z}$--module. In this case  for any proper subgroup  $H$ such that  $L \leq H$ 
it is existed  a series of $\mathbb{Z}$$G$--submodules  
$$
 0 \hspace{0,2cm} \leq \hspace{0,2cm} C_{A}(H) \hspace{0,2cm} \leq \hspace{0,2cm}
A_{1} \hspace{0,2cm} \leq \hspace{0,2cm} A
$$
 such that the additive group of  $A_{1}/ C_{A}(H)$ is either an abelian Chernikov group or trivial,  the additive group of
$A/A_{1}$ is an abelian torsion-free group of finite  $0$--rank. Therefore we can construct the series of 
$\mathbb{Z}$$G$--submodules 
$$
0 \hspace{0,2cm} \leq \hspace{0,2cm} C_{A}(H) \hspace{0,2cm} \leq \hspace{0,2cm}
A_{1} \hspace{0,2cm} \leq \hspace{0,2cm} A_{2} \hspace{0,2cm} \leq \hspace{0,2cm} \cdots \hspace{0,2cm} \leq \hspace{0,2cm} A_{n-1}
\hspace{0,2cm} \leq \hspace{0,2cm} A_{n} = A,
$$
such that  the additive group of  $A_{1}/ C_{A}(H)$ is either abelian Chernikov or trivial, 
 $A_{k+1}/A_{k}$, $k = 1, \cdots, n-1$,  are  $G$--rationally irreducible and  the additive groups of  $A_{k+1}/A_{k}$ are abelian torsion--free groups of finite  $0$--rank.

 Now we consider  the series of  $\mathbb{Z}$$G$--submodules 
$$
\overline{ 0} = C_{A}(H)/C_{A}(H) \leq  
A_{1}/C_{A}(H)  \leq  A_{2}/C_{A}(H) \leq  \cdots  
$$
$$
\leq  A_{n-1}/C_{A}(H)  \leq  A_{n}/C_{A}(H) = A/C_{A}(H).
$$
If the  additive group of $A_{1}/C_{A}(H)$ is non--trivial then  $A_{1}/C_{A}(H)$ is the  union of finite characteristic subgroups $B_{m}/C_{A}(H)$, $ m=1, 2, \cdots$, 
and for each  $ m=1, 2, \cdots$, we have that   $G/C_{G}(B_{m}/C_{A}(H))$ is finite. By lemma 2.2  $G$ has not proper subgroups of finite index.  Then  $G = C_{G}(B_{m}/C_{A}(H))$ 
for each  $ m=1, 2, \cdots$. Therefore  $[G, B_{m}] \leq C_{A}(H)$ for each  $ m=1, 2, \cdots$. It follows that  
$[G,A_{1}] \leq C_{A}(H)$. So $G$ acts trivially at the factor $A_{1}/C_{A}(H)$. 
For any  $k = 1, \cdots, n-1$,  we have the isomorphism  of  $\mathbb{Z}$$G$--modules
$$
(A_{k+1}/C_{A}(H))/(A_{k}/C_{A}(H)) \simeq A_{k+1}/A_{k}.
$$
 The quotient group  $G/C_{G}(A_{k+1}/A_{k}), k = 1, \cdots, n-1,$
 can be considered as
an  irreducible subgroup of $GL_{r}(\mathbb{Q})$.  By corollary 3.8  \cite{BAW} $G/C_{G}(A_{k+1}/A_{k})$ is a soluble group for each  $k = 1, \cdots, n-1$. By A.I.Maltzev theorem  (lemma 3.5 \cite{BAW}) $G/C_{G}(A_{k+1}/A_{k})$, $k = 1, \cdots, n-1$, are almost abelian. 
Since  $G$ has not proper subgroups of finite index then $G/C_{G}(A_{k+1}/A_{k})$, $k = 1, \cdots, n-1$, are abelian. From  $D \leq C_{G}(A_{k+1}/A_{k})$ it follows that for any  $k = 1, \cdots, n-1$, 
the quotient group  $G/C_{G}(A_{k+1}/A_{k})$ is either  isomorphic to a quasi--cyclic $q$--group for some prime $q$ or trivial.
For each  $k = 1, \cdots, n-1$,  the additive group of  $A_{k+1}/A_{k}$ is an abelian  $A_{4}$--group \cite {MAL}. Since a periodic subgroup of automorphisms group of a soluble  $A_{4}$--group is 
finite \cite{KSM}, then $G/C_{G}(A_{k+1}/A_{k})$, $k = 1, \cdots, n-1$, are trivial. 
Therefore  $G$ acts trivially at every factor of the series 
$$
\overline{0}  \leq  
A_{1}/C_{A}(H)  \leq  A_{2}/C_{A}(H) \leq  \cdots  \leq  A_{n-1}/C_{A}(H)
 \leq A_{n}/C_{A}(H) = A/C_{A}(H).
$$

 At first  we consider the case where there exists a subgroup  $M$ such that $D \leq M$,    
$M/D$ is a finite subgroup of  $G/D$, $ A/C_{A}(M)$ is not an artinian $\mathbb{Z}$--module and for any proper subgroup $H$ of $G$ such that   $M \leq H$  
the equalities  $A_{1} = C_{A}(H)$, $A_{2} = A$ are valid. As we proved   $G$ acts trivially at the factor  $A_{2}/A_{1} = A/C_{A}(H)$. Therefore   $[G,A] \leq C_{A}(H)$. From the choice of  $H$ it follows that $[G,A] \leq C_{A}(G)$. 
So $G$ acts trivially at every factor of the series 
$0 \hspace{0,2cm} \leq \hspace{0,2cm} C_{A}(G) \hspace{0,2cm} \leq \hspace{0,2cm}A.$
 It follows that  $G$ is an abelian group. Therefore 
$G$ is isomorphic to a quasi--cyclic $q$--group. Otherwise we can choose a subgroup $H$ such that  $H/D$ is a finite subgroup of  $G/D$, $ A/C_{A}(H)$ 
is not an artinian  $\mathbb{Z}$--module and even if one from the  equalities $A_{1} = C_{A}(H)$, $A_{2} = A$ is not valid. Then  by Kaluzhnin theorem (p. 144 \cite {KAM})
 $G$ is a nilpotent group of step  $\leq n-1$.
Let  
$$
\langle 1 \rangle = Z_{0} \hspace{0,2cm} \leq \hspace{0,2cm} Z_{1} \hspace{0,2cm} \leq \hspace{0,2cm} Z_{2} \hspace{0,2cm} \leq 
\hspace{0,2cm} \cdots \hspace{0,2cm} \leq \hspace{0,2cm} Z_{l-1} \hspace{0,2cm} \leq \hspace{0,2cm} Z_{l} = G
$$
 be an upper central series of  $G$, where  $l \leq n-1$. 
From $D \leq Z_{l-1}$ it follows that  $G/Z_{l-1}$ is  isomorphic to a quasi--cyclic $q$--group for some prime $q$. Therefore  
$G/Z_{l-2}$ is an extension of a central subgroup by a  quasi--cyclic $q$--group.  Then   $G/Z_{l-2}$ is abelian. By lemma 2.2   $G/Z_{l-2}$ is  isomorphic to a quasi--cyclic $q$--group for some prime $q$. 
If we continue  similarly at  the step with the number   $l-1$ we obtain that  $G/Z_{1}$ 
is  isomorphic to a quasi--cyclic $q$--group for some prime $q$.  It follows that  $G$ is an extension of a central subgroup by a  quasi--cyclic $q$--group.
Therefore  $G$ is an abelian group and $G$ is  isomorphic to a quasi--cyclic $q$--group $C_{q^{\infty}}$ for some prime $q$.

Now we consider the case where for every  proper subgroup  $H$ such that  $D \leq H$  the quotient module   $ A/C_{A}(H)$ is  an artinian $\mathbb{Z}$--module.
It follows from theorem 3.1 \cite {OYD} that  $G$ is  isomorphic to a quasi--cyclic $q$--group $C_{q^{\infty}}$ for some prime $q$. The theorem is proved.

It should be noted that theorem 2.2 is the special  case of the main theorem of  \cite{KSC}. 

{\bf Theorem 2.3.} 
 Let $A$ be an  $ \bf R$$G$--module where  $G$ is an infinite group, $G \not = G'$, $ \bf R$ is an associative ring. 
If $A/C_{A}(G)$ is an infinite $ \bf R$--module and $A/C_{A}(H)$ is a finite $ \bf R$--module for every proper subgroup $H$ of $G$ then $G$ is isomorphic to a quasi--cyclic $q$--group  $C_{q^{\infty}}$ for some prime $q$.

{\bf Proof.}
Let  $D$ be a derived subgroup of $G$. By lemma 2.2    $G/D$ is isomorphic to a quasi--cyclic $q$--group for some prime $q$. 
Let $H/D$ be any finite subgroup of $G/D$.
Since   $H$  is a proper subgroup of  $G$ then $A/C_{A}(H)$ is finite.  Therefore   
 $G/C_{G}(A/C_{A}(H))$ is finite.  As  by lemma 2.2 $G$ has not proper subgroups of finite index then $G = C_{G}(A/C_{A}(H))$. It follows that $[G, A] \leq C_{A}(H)$. 
From the choice of  $H$  it follows that  $[G,A]  \leq  C_{A}(G)$. So   
  $G$ acts trivially at every factor of the series 
$0  \leq  C_{A}(G) \leq A$.  By Kaluzhnin theorem (p. 144 \cite {KAM})  
  $G$ is an abelian group. Therefore $G$   is  isomorphic to a quasi--cyclic $q$--group $C_{q^{\infty}}$ for some prime $q$.    
The theorem is proved.

In \cite{KSC} the authors have constructed the example of a module with the prescribed conditions.

\end{document}